\documentclass[11pt]{articlefederico}
\usepackage{amsthm}
\usepackage{amsmath}
\usepackage{amssymb}

\usepackage{graphicx}

\newcommand{\B}{{\cal B}}
\newcommand{\C}{\mathbb{C}}

\newcommand{\R}{\mathbb{R}}

\begin{document}
\title{\textsf{Transversal and cotransversal matroids via the Lindstr\"om lemma.}}
\author{\textsf{Federico Ardila\footnote{\textsf{Department of Mathematics, San
Francisco State University, federico@math.sfsu.edu. \newline
\hspace{2cm} Supported by NSF grant DMS-9983797. }}}}
\date{}
\maketitle

\begin{abstract}
It is known that the duals of transversal matroids are precisely the
strict gammoids. The purpose of this short note is to show how the
Lindstr\"om-Gessel-Viennot lemma gives a simple proof of this
result.
\end{abstract}
%
%
%
%

\bigskip
\bigskip
\noindent{\bf \textsf{\Large{1}}}
\bigskip

\noindent{\bf \textsf{Matroids and duality.}} A \emph{matroid}
$M=(E,\B)$ is a finite set $E$, together with a non-empty collection
$\B$ of subsets of $E$, called the \emph{bases} of $M$, which
satisfy the following axiom: If $B_1, B_2$ are bases and $e$ is in
$B_1-B_2$, there exists $f$ in $B_2-B_1$ such that $B_1-e\cup f$ is
a basis.

If $M=(E,\B)$ is a matroid, then $\B^*=\{E-B \, | \, B \in \B\}$
is also the collection of bases of a matroid $M^*=(E,\B^*)$,
called the \emph{dual} of $M$.

\medskip

\noindent{\bf \textsf{Representable matroids.}} Matroids can be
thought of as a combinatorial abstraction of linear independence.
If $V$ is a set of vectors in $\R^n$ and $\B$ is the collection of
maximal linearly independent sets of $V$, then $M=(V,\B)$ is a
matroid. Such a matroid is called \emph{representable} over $\R$,
and $V$ is called a \emph{representation} of $M$.


\medskip

\noindent{\bf \textsf{Transversal matroids.}} Let $A_1, \ldots,
A_r$ be subsets of $[n]=\{1,\ldots,n\}$. A \emph{transversal}
(also known as \emph{system of distinct representatives}) of
$(A_1, \ldots, A_r)$ is a subset $\{e_1, \ldots, e_r\}$ of $[n]$
such that $e_i$ is in $A_i$ for each $i$. The transversals of
$(A_1, \ldots, A_r)$ are the bases of a matroid on $[n]$. Such a
matroid is called a \emph{transversal matroid}, and $(A_1, \ldots,
A_r)$ is called a \emph{presentation} of the matroid. This
presentation can be encoded in the bipartite graph $H$ with
``left" vertex set $L=[n]$, ``right" vertex set
$R=\{\widehat{1},\ldots,\widehat{r}\}$, and an edge joining $j$
and $\widehat{i}$ whenever $j$ is in $A_i$. The transversals are
the $r$-sets in $L$ which can be matched to $R$. We will denote
this transversal matroid by $M[H]$.

\medskip

\noindent{\bf \textsf{Strict gammoids.}} Let $G$ be a directed
graph with vertex set $[n]$, and let $A=\{v_1, \ldots, v_r\}$ be a
subset of $[n]$. We say that an $r$-subset $B$ of $[n]$ \emph{can
be linked to $A$} if there exist $r$ vertex-disjoint directed
paths whose initial vertex is in $B$ and whose final vertex is in
$A$. We will call these $r$ paths a \emph{routing} from $B$ to
$A$. The collection of $r$-subsets which can be linked to $A$ are
the bases of a matroid denoted $L(G,A)$. Such a matroid is called
a \emph{strict gammoid}.

We can assume that the vertices in $A$ are sinks of $G$;
\emph{i.e.}, that there are no edges coming out of them. This is
because the removal of those edges does not affect the matroid
$L(G,A)$.

%

\bigskip
\bigskip
\noindent{\bf \textsf{\Large{2}}}
\bigskip

\noindent{\bf \textsf{Representations of transversal matroids.}}
Consider a collection of algebraically independent $\alpha_{ij}$s
for $1 \leq i \leq r, 1 \leq j \leq n$. Let $M$ be a transversal
matroid on the set $[n]$ with presentation $(A_1, \ldots, A_r)$. Let
$X$ be the $r \times n$ matrix whose $(i,j)$ entry is $-\alpha_{ij}$
if $j \in A_i$ and $0$ otherwise. The columns of $X$ are a
representation of $M$.

To see this, consider the columns $j_1, \ldots, j_r$. They are
independent when their determinant is non-zero. As soon as one of
the $r!$ summands in the determinant is non-zero, the determinant
itself will be non-zero, by the algebraic independence of the
$\alpha_{ij}$s. But the summand $\pm X_{\sigma_1j_1}\cdots
X_{\sigma_rj_r}$ (where $\sigma$ is a permutation of $[r]$) is
non-zero if and only if $j_1 \in A_{\sigma_1}, \ldots, j_r \in
A_{\sigma_r}$. So the determinant is non-zero if and only if $\{j_1,
\ldots, j_r\}$ is a transversal. The desired result follows.

We will find it convenient to choose a transversal $j_1 \in A_1,
\ldots, j_r \in A_r$ ahead of time, and normalize the rows to have
$-\alpha_{ij_i}=1$ for $1 \leq i \leq r$.

\medskip

\noindent\textit{Example 1}. Let $n=6$ and $A_1 = \{1,2,3\}, A_2 =
\{2,4,5\}, A_3 = \{3,5,6\}.$  The corresponding bipartite graph
$H$ is shown below.

\begin{figure}[h]
\centering
\includegraphics[height=3.5cm]{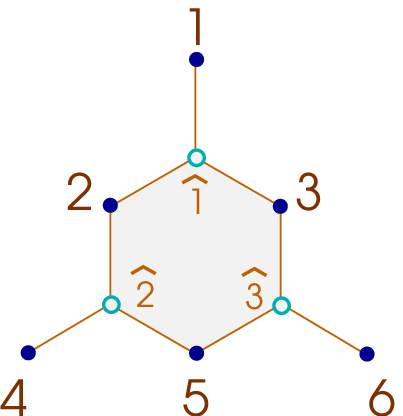}
\label{fig:h}
\end{figure}

If we choose the transversal $1 \in A_1, 2 \in A_2, 3 \in A_3$, we
obtain a representation for the transversal matroid $M[H]$, given
by the columns of the following matrix:
\[
X=\left( \begin{array}{cccccc}
1 & -a & -b & 0 & 0 & 0 \\
0 & 1 & 0 & -c & -d & 0 \\
0 & 0 & 1 & 0 & -e & -f
\end{array}
\right)
\]

\medskip

\noindent{\bf \textsf{Representations of strict gammoids.}} Let
$M=L(G,A)$ be a strict gammoid. Say $G$ has vertex set $\{1,
\ldots, n\}$ and $A=\{a_1, \ldots, a_{n-r}\}$. Assign
algebraically independent weights smaller than $1$ to the edges of
$G_n$. For $1 \leq i \leq n-r$ and $1 \leq j \leq n$, let $p_{ij}$
be the sum of the weights of all finite paths\footnote{The weight
of a path is defined to be the product of the weights of its
edges. The sum converges since the weights are less than $1$.}
from vertex $i$ to vertex $j$. Let $Y$ be the $(n-r) \times n$
matrix whose $(i,j)$ entry is $p_{ji}$. The columns of $Y$ are a
representation of $M$.

This is a direct consequence of the Lindstr\"om lemma or
Gessel-Viennot method, which tells us that the determinant of the
matrix with columns $j_1, \ldots, j_{n-r}$ is equal to the signed
sum\footnote{The sign is determined by the permutation that matches
the starting and ending points of the paths in the routing.} of the
routings from $\{j_1, \ldots, j_{n-r}\}$ to
$\{a_1,\ldots,a_{n-r}\}$. This signed sum is non-zero if and only if
it is non-empty.

\medskip

\noindent\textit{Example 2}. Consider the graph $G$ shown below,
where all edges point down, and the set of sinks $A=\{4,5,6\}$.

\begin{figure}[h]
\centering
\includegraphics[height=3.5cm]{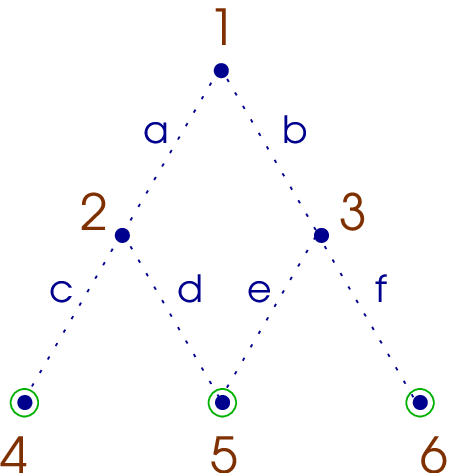}
\label{fig:g}
\end{figure}

The representation we obtain for the strict gammoid $L(G,A)$ is
given by the columns of the following matrix:
\[
Y = \left(
\begin{array}{cccccc}
ac & c & 0 & 1 & 0 & 0 \\
ad+be & d & e & 0 & 1 & 0 \\
bf & 0 & f & 0 & 0 & 1
\end{array}
\right)
\]

Notice that the rowspaces of $X$ and $Y$ are orthogonally
complementary in $\R^6$. That is, essentially, the punchline of
this story.

\bigskip
\bigskip
\noindent{\bf \textsf{\Large{3}}}
\bigskip

\noindent{\bf \textsf{Representations of dual matroids.}} If a rank
$r$ matroid $M$ is represented by the columns of an $r \times n$
matrix $A$, we can think of $M$ as being represented by the
$r$-dimensional subspace $V = \textrm{rowspace}(A)$ in $\R^n$. The
reason is that, if we consider any other $r \times n$ matrix $A'$
with $V = \textrm{rowspace}(A')$, the columns of $A'$ also represent
$M$.

This point of view is very amenable to matroid duality. If $M$ is
represented by the $r$-dimensional subspace $V$ of $\R^n$, then
the dual matroid $M^*$ is represented by the $(n-r)$-dimensional
orthogonal complement $V^*$ of $\R^n$.

\medskip

\noindent{\bf \textsf{Digraphs with sinks and bipartite graphs
with complete matchings.}}  From a directed graph $G$ on the set
$[n]$ and a set of $n-r$ sinks $A \subseteq [n]$ of $G$, we can
construct a bipartite graph $H$ as follows. The left vertex set is
$[n]$, and the right vertex set is a copy $[\widehat{n}] -
\widehat{A}$ of $[n]-A$. We join $\widehat{u}$ and $u$ for each $u
\in [n]-A$, and we join $\widehat{u}$ and $v$ whenever $u
\rightarrow v$ is an edge of $G$. This graph $H$ has the obvious
complete matching between $\widehat{u}$ and $u$. Conversely, if we
are given the bipartite graph $H$ with a complete matching, it is
clear how to recover $G$ and $A$.

Observe that if we start with the directed graph $G$ and sinks $A$
of Example 1, we obtain the bipartite graph $H$ of Example 2.
%
%

\medskip

\noindent{\bf \textsf{Duality of transversal matroids and strict
gammoids.}} Now we show that, in the above correspondence between a
graph $G$ with sinks $A$ and a bipartite graph $H$ with a complete
matching, the strict gammoid $L(G,A)$ is dual to the transversal
matroid $M[H]$. We have constructed a subspace of $\R^n$
representing each one of them, and now we will see that they are
orthogonally complementary, as observed in Examples 1 and 2.

%
%
%
%

Our representation of $M[H]$ is given by the columns of the $r
\times n$ matrix $X$ whose $(i,i)$ entry is $1$, and whose $(i,j)$
entry is $-\alpha_{ij}$ if $i \rightarrow j$ is an edge of $G$ and
$0$ otherwise. Think of the $\alpha_{ij}$s as weights on the edges
of $G$. A vector $y \in \C^n$ is in the $(n-r)$-dimensional null
space of $X$ when, for each vertex $i$ of $G$,
\begin{equation}\label{eqn}
y_i = \sum_{j \in N(i)} \alpha_{ij}y_j.
\end{equation}
Here $N(i)$ denotes the set of vertices $j$ such that $i \rightarrow
j$ is an edge of $G$.

As before, let $p_{ia}$ be the sum of the weights of the finite
paths from $i$ to $a$ in $G$. Our representation $Y$ of $L(G,A)$ has
rows $(y_1, \ldots, y_n) = (p_{1a}, \ldots, p_{na})$ (for $a \in
A$). Clearly, each row of $Y$ is a solution to $(1)$, so
$\textrm{rowspace}(Y) \subseteq \textrm{nullspace}(X)$. But these
two subspaces are $(n-r)$-dimensional, so they must be equal, as we
wished to show.
%
%
This completes our proof of the theorem that the strict gammoids are
precisely the cotransversal matroids.

\bigskip
\bigskip
\noindent{\bf \textsf{\Large{4}}}
\bigskip


For more information on matroid theory, Oxley's book \cite{Oxley} is
a wonderful place to start. The representation of transversal
matroids shown here is due to Mirsky and Perfect \cite{Mirsky}. The
representation of strict gammoids that we use was constructed by
Mason \cite{Mason} and further explained by Lindstr\"om
\cite{Lindstrom}\footnote{It is in this context that he discovered
what is now known as the Lindstr\"om lemma or Gessel-Viennot method
\cite{Gessel}. This method was also used earlier by Karlin and
MacGregor \cite{Karlin}.}. The theorem that strict gammoids are
precisely the cotransversal matroids is due to Ingleton and Piff
\cite{Ingleton}. Our proof of this result appears to be new.

%

%

This note is a small side project of \cite{Ardila}. While studying
the geometry of flag arrangements and its implications on the
Schubert calculus, we were led to study a specific family of strict
gammoids which starts with Example 2. I would like to thank Sara
Billey for several helpful discussions, and Laci Lovasz and Jim
Oxley for help with the references.

\small{

}

\end{document}